\documentclass[12pt]{amsart}

\usepackage{amsthm}

\usepackage{epsf}

\usepackage{amssymb}

\usepackage{latexsym}

\usepackage{amsmath}

\newtheorem{theorem}{Theorem} 

\newtheorem{definition}[theorem]{Definition}

\newtheorem{lemma}[theorem]{Lemma}

\newtheorem{proposition}[theorem]{Proposition}

\newtheorem{remark}[theorem]{Remark}

\newtheorem{notation}[theorem]{Notation}

\def\C{{\mathbb C}}

\def\T{{\mathbb T}}

\def\Z{{\mathbb Z}}

\def\P{{\mathbb P}}

\def \ta{\tau}

\def \ta1{\tau_1}

\def \g{\gamma}

\def \G{\Gamma}

\def \ri{\rightarrow}

\def \CP{\mathbb{C} \mathbb{P}}

\def\isom{{\cong}}

\newcommand\begintable[1][] {{}}

\long\def\forget#1\forgotten{}

\newif\ifXY 

\XYtrue     

\def\CP{\mathbb{CP}}

\def\Z{\mathbb{Z}}

\begin{document}

\renewcommand{\subjclassname}{%

       \textup{2000} Mathematics Subject Classification}

\date{\today}

\title[fundamental group of the complement of Cayley's singularities]{The fundamental group of the complement of the Cayley's singularities$^{1}$}

\author[Amram, Dettweiler, Friedman, Teicher]{Amram Meirav, Dettweiler
 Michael, Friedman Michael  and  Teicher Mina}

\stepcounter{footnote}

\footnotetext{This work was supported in part by the Edmund Landau
Center for Research in Mathematics, Hebrew University; the Emmy
Noether Research Institute for Mathematics, Bar-Ilan University; the
Minerva Foundation (Germany), the EU-network
HPRN-CT-2009-00099(EAGER), and the Israel Science Foundation grant
$\#$ 8008/02-3 (Excellency Center ``Group Theoretic Methods in the
Study of Algebraic Varieties").}

\address{Meirav Amram, Einstein Mathematics Institute, Hebrew University,
Jerusalem, Israel; Department of Mathematics, Bar-Ilan University,
Israel} \email{ameirav@math.huji.ac.il, meirav@macs.biu.ac.il}

\address{Michael Dettweiler, IWR, Heidelberg University, Germany}

\email{michael.dettweiler@iwr.uni-heidelberg.de}

\address{Michael Friedman, Mina Teicher, Department of Mathematics, Bar-Ilan University, Israel}

\email{fridmam@macs.biu.ac.il, teicher@macs.biu.ac.il}

\begin{abstract}
Given a singular surface $X$, one can extract information on it by investigating
 the fundamental group $\pi_1(X - Sing_X)$. However, calculation of this group is non-trivial, but it can
 be simplified if a certain invariant of the branch curve of $X$ -- called the braid monodromy factorization -- is known.
 This paper shows, taking the Cayley cubic as an example, how this fundamental group can be computed
 by using braid monodromy techniques (\cite{16}) and their liftings. This is
one of the first examples that uses these techniques to calculate this sort of
fundamental group.
\end{abstract}

\keywords{Singularities, Coverings, Fundamental Groups, Surfaces,
Mapping Class Group \\ MSC Classification: 14B05, 14E20, 14H30,
14Q10}

\maketitle

\section{On the Cayley cubic}\label{overview}

The classification of singular cubic surfaces in $\C\P^3$ was done
in the 1860's, by Schl\"{a}fli \cite{Sch} and Cayley \cite{Cay}.
Surface XVI in Cayley's classification is now called the Cayley
cubic, and when embedded in $\C\P^3$, it is defined by the following
equation
$$
4 (X^3 + Y^3 + Z^3 + W^3) - (X + Y + Z + W)^3=0.
$$
It has four singularities, which are ordinary double points. Cayley
noticed that this surface is the unique cubic surface having four
ordinary double points, which is the maximal possible number of
double points for a cubic surface (see, for example, Salmon's book
\cite{Sal}). Note that the Cayley cubic is invariant under the symmetric group
$Sym_4$ and it contains exactly nine lines, six of which connect the four nodes pairwise and the
other three of which are coplanar (see, for example \cite[Section 4.1.3]{hunt}).

Let us denote this surface by $C$ and by $Sing_C = Sing$ the set of
the four nodes. We are interested in the fundamental group of the
complement of the set of singularities in the Cayley cubic
$\pi_1(C-Sing)$. A direct computation of this group is elementary.
Consider the smooth del-Pezzo surface $S_6$ of degree $6$. The
Cremona involution $Cr$ is the regular automorphism of $S_6$ and has
$4$ fixed points, which we denote as $\text{Fix}(Cr)$. The factor
$S_6/Cr$ is the Cayley cubic; singular points of the Cayley cubic
are images of fixed points of $Cr$ on $S_6$. Therefore, the
universal cover of ($C - Sing$) is ($S_6 - \text{Fix}(Cr)$), and
hence the fundamental group $\pi_1(C-Sing)$ is
$\mathbb{Z}/2\mathbb{Z}$.

However, for a general singular surface $X$ in $\C\P^3$, there is no
general method for computing the fundamental group $\pi_1(X -
Sing_X)$. We present here two other different approaches for this
problem, demonstrating them on the Cayley cubic. We compute first
the braid monodromy factorization of the branch curve of the Cayley
cubic in $\C\P^2$, based on the braid monodromy techniques of
Moishezon-Teicher (\cite{16}, \cite{17}). We then apply two methods
in order to compute this fundamental group. The first method
consists of lifting the factorization to a factorization in the
mapping class group, from which we can find the desired group. The
second method is based on {\cite{Manf}} and finds the fundamental
group using the Reidemeister-Schreier method (\cite{MKS}).

The paper is divided as follows.
  In Section \ref{sec:5} we compute the braid monodromy factorization of the
branch curve of the  Cayley cubic $C$ and the fundamental group of
the complement of the branch curve. In Section \ref{sec:6} the
fundamental group of the Cayley surface minus the singular points
is computed using the results from Section \ref{sec:5}.

\section{The factorization $\Delta^2_6$ and the fundamental group $\pi_1(\C\P^2 - \bar{S})$}\label{sec:5}

In this section we give the braid monodromy factorization of the
branch curve $S$ in $\C^2$. We also present the fundamental groups
$\pi_1(\C^2 - {S})$ and $\pi_1(\C\P^2 - \bar{S})$.

We begin with a few basic notations. Let $\bar{S}$ be a branch curve
of a surface $C$ in $\C\P^2$, and let $l_{\infty}$ be a line in
$\C\P^2$, transverse to $\bar{S}$. Let $S = \bar{S} -
\bar{S} \cap l_{\infty}$. Take a projection $\pi : \C\P^2
\rightarrow \C\P^1$, and let $\pi_{aff}: \C^2 \rightarrow \C$ be the corresponding
affine generic projection. Then there is a finite subset of points
$Z \subset \C$, which is the projection on $\C$ of the nodes and
cusps of $S$ and the branch points of $\pi_{aff} |_{S}$. Above each
point of $Z \subset \C$ there is just one singular point of
$\pi_{aff} |_{S}$. Let $\{ \delta_{i} \}$ be a basis of
non-intersecting loops in $\C - Z$ around each point of $Z$,
starting from $x_{0} \in \C - Z$.

Consider now a  closed disc $D \subset \mathbb{R}^2,\,K$ a finite set of points in it, and look at $\mathcal{B}$
the group of all the diffeomorphisms $\beta$ such that $\beta(K) = K$ and $\beta|_{\partial D} = id$. We say
that $\beta_1, \beta_2 \in \mathcal{B} $ are equivalent if they induce the same automorphism of $\pi_1(D \setminus K, u).$ The quotient of $\mathcal{B}$ by this equivalence relation is called the braid
group $B_n = B_n[D,K]\,\, (n = \#K)$.

Now we use the theorem of Zariski \cite{Z}: {\em Let $z \in Z$ and
$\delta$ be a loop in $\C - Z$ around $z$. Then there is a
\emph{braid monodromy} action $\varphi: \pi_{1}(\C - Z, x_{0})
\rightarrow B_{n}[\C_{x_{0}}, \C_{x_{0}} \cap S]$, s.t. $\C
_{x_{0}}$ is the fiber of $\pi_{aff}$ over $x_{0}$ and $B_n$ is the
braid group}.

Assume that $C_{x_0}\cap S$ is on the $x$-axis. Then we define a halftwist $Z_{ij}$ as the exchange of the positions of two points $i$ and $j$, which occurs as follows: we take a tubular
neighborhood of a path which connects $i$ and $j$ below the
$x$-axis, then we let $i$ and $j$ rotate in a counterclockwise
manner along the boundary of this neighborhood, until they
exchange their position.
\begin{notation}
We denote by $Z_{ij}$ (resp. $\bar{Z}_{ij}$) the counterclockwise
halftwist of $i$ and $j$ below (resp. above) the axis. $Z_{ij}^2$
is a fulltwist of $i$ around $j$. $Z^{2}_{i,jj'}=Z^{2}_{ij'}
Z^{2}_{ij}$ is the fulltwist of $i$ around $j$ and $j'$. In a
similar way, we define also $Z^2_{ii',jj'}=Z^2_{i',jj'}
Z^2_{i,jj'}$  and $Z^{3}_{i',jj'}=Z^{3}_{i'j } Z^{3}_{i'j'}
{(Z^{3}_{i'j})}^{Z^2_{jj'}}$.
\end{notation}

In the case when the singular point above $z$ is a branch point of
$\pi_{aff}$, a node, a cusp of $S$, or a point of tangency of two
branches of the curve, then $\varphi{(\delta)} = {H}^{\epsilon}$,
where $H$ is a halftwist and $\epsilon = 1, 2, 3, 4$
(respectively).

When $z$ can be given locally as an intersection on $m$ lines,
then $\varphi{(\delta)} = \Delta^2_m$, when $\Delta^2_m$ is a
$360$ degree rotation of the $m$ points in the fiber.

More details for explicit computations and technical methods
appear in \cite{16} and \cite{17}.

\begin{definition}
The braid monodromy (=BMF) w.r.t. $S,\pi,u$ is the following
factorization
$$
\Delta^2_S = \prod_i \varphi(\delta_i).
$$
\end{definition}

\begin{remark} \label{remArt} Let $\Delta^2$ be the generator of the center of the braid group $B_n[D,K]$ .
Then, by a theorem of Artin (see \cite{16}), $\Delta^2 =
\prod\varphi(\delta_i).$ Note that $\Delta^2$ is a $360$ degree
rotation of the disc $D$.
\end{remark}

\subsection{BMF of the Cayley cubic}

Denote the Cayley cubic as $C$ and the set of the four nodes as
$Sing$. We aim to compute the fundamental group $\pi_{1}(C -
Sing)$. Denote also $V(d,c,n)$ the variety of degree $d$ plane curves with $c$ cusps
and $n$ nodes. It can be seen easily that $S$ -- the branch curve of $C$ -- is in $V(6,6,4)$ and has
$4$ branch points with respect to a generic projection to $\C\P^1$. Note also that the dual curve $S^{\vee}$ belongs
to $V(4,0,3)$. Since the variety $V(4,0,3)$ is irreducible, $V(6,6,4)$ is also irreducible. Therefore, we can pick
any curve $S' \in V(6,6,4)$, since $S$ and $S'$ would be isotopic and thus their braid monodromy factorizations
would be equivalent (see \cite{KuTe}).  In particular, the fundamental groups of the complement of $S$ and $S'$ will
be isomorphic. In fact, we do not find explicitly a curve $S' \in V(6,6,4)$. We will start from a more basic branch curve of a degenerated surface and ``regenerate" it; this process will recover the braid monodromy factorization
of a curve $S' \in V(6,6,4)$.

Let us consider a union of three planes meeting at a point: we call this surface the degenerated surface.
The branch curve $S_0$ of this surface is an arrangement of three lines meeting at one point, one
of which is set to be the ``diagonal line". Denote by $U \subset \C^2$ a small neighborhood
of the singular point.

%


The next step is to apply the regeneration process. When
regenerating a singular configuration consisting of lines and
conics, the final stage in the regeneration process involves
doubling each line. Let $x_0 \in \mathbb{C}$ be a generic point such that $\pi_{aff}^{-1}(x_0) \cap S_0 \in U \setminus Sing\,S_0$. In the regeneration process each point of $K_0 = \pi_{aff}^{-1}(x_0) \cap S_0$ corresponding to a
line labelled $i$ is replaced by a pair of points, labelled $i$ and
$i'$. The purpose of the regeneration rules is to explain how the
braid monodromy behaves when lines are doubled in this manner. We
denote $H(z_{i,j})$ by $Z_{i,j}$ (where $z_{i,j}$ is a path
connecting points in $K$).

The rules are (see \cite[pp. 336-337]{19}):
\begin{enumerate}
\item \textbf{First regeneration rule}: The regeneration of a
branch point of any conic -- any branch point regenerates into two branch points:\\
A factor of the braid monodromy of the form $Z_{i,j}$ is replaced
in the regeneration by $Z_{i',j}\cdot
\overset{(j)}{\underline{Z}}_{i,j'}$\medskip
\item \textbf{Second regeneration rule}: The regeneration of a node -- any node regenerates into two (or four) nodes:\\
A factor of the form $Z^2_{ij}$ is replaced by a factorized
expression $Z^2_{ii',j} := Z^2_{i'j}\cdot Z^2_{ij}$ ,\\
$Z^2_{i,jj'} := Z^2_{ij'}\cdot Z^2_{ij}$ or by $Z^2_{ii',jj'} :=
Z^2_{i'j'}\cdot Z^2_{ij'}Z^2_{i'j}\cdot Z^2_{ij}$. \medskip \item
\textbf{Third regeneration rule}: The regeneration of a tangent
point -- any tangent point regenerates into three cusps:\\
A factor of the form $Z^4_{ij}$ in the braid monodromy factorized
expression is replaced by\\ $Z^3_{i,jj'} :=
(Z^3_{ij})^{Z_{jj'}}\cdot (Z^3_{ij}) \cdot
(Z^3_{ij})^{Z^{-1}_{jj'}}$.
\end{enumerate}

\begin{figure}

\begin{minipage}{\textwidth}

\begin{center}

\epsfbox{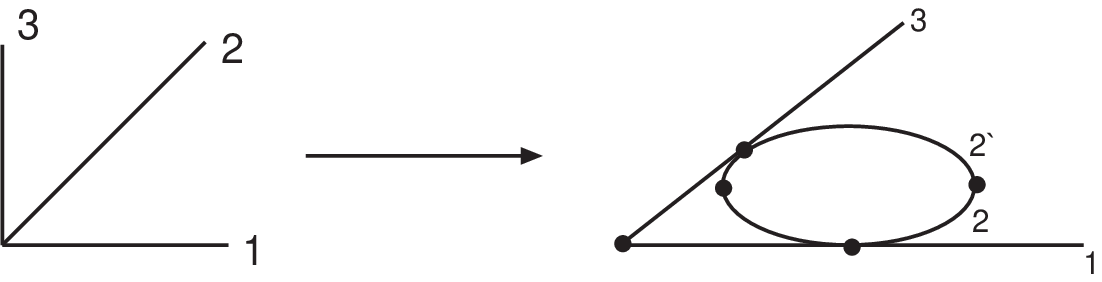}

\end{center}

\end{minipage}

\caption{}\label{conic}

\end{figure}

The initial braid monodromy factorization of the branch curve (which is a union of three
lines meeting at a point) is $\Delta^2_3$. We first regenerate  the ``diagonal line"
to a smooth conic which is tangent to the two other lines (see
\cite[Lemma 1]{19}), as depicted in Figure \ref{conic}. The braid
monodromy factorization of this arrangement is
$$
(Z^2_{1,3})^{Z^{-1}_{2',3}}\cdot
Z_{2,2'}^{Z^{-2}_{1,2}\bar{Z}^2_{2,3}} \cdot
(Z_{2,3}^4)^{Z^2_{2,2'}}\cdot Z_{1,2'}^4 \cdot Z_{2,2'}.
$$

\begin{figure}




\begin{center}

\epsfbox{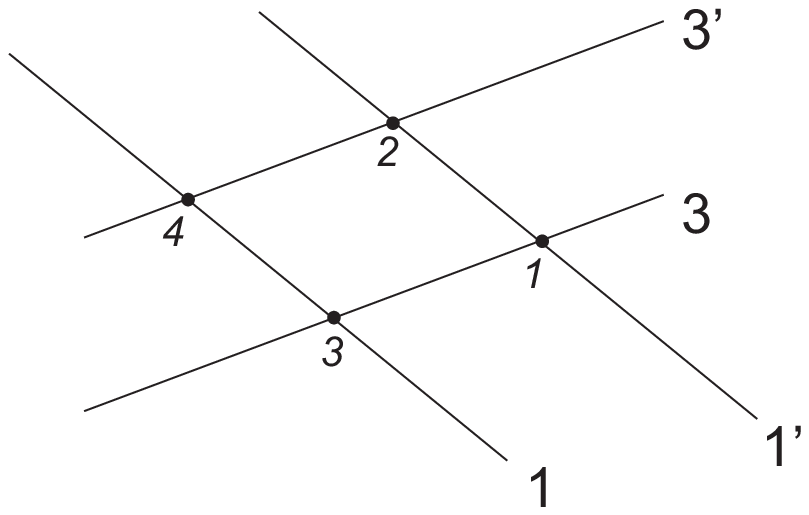}

\end{center}


\caption{}\label{4lines}

\end{figure}

We now regenerate the remaining two lines. By the second
regeneration rule, the node is regenerated into four nodes, and
each tangency point regenerates into three cusps. For example, the
regeneration in a neighborhood of the node is depicted in Figure
\ref{4lines}.

 We end up with a curve $\tilde{S}$ which has -- in the (regenerated) neighborhood $U$ --  $6$ cusps, $4$ nodes and $2$ branch points.
 The resulting factorization is
 \begin{eqnarray} \label{del1}
\widetilde{\Delta} = & {(Z_{1' \; 3}^2)}^{Z^{-2}_{2', 3 \; 3'}}
\cdot (Z_{1' \; 3'}^2)^{Z_{1' \; 3}^2 Z^{-2}_{2', 3 \; 3'}} \cdot
{(Z_{1 \; 3}^2)}^{Z^{-2}_{2', 3 \; 3'}} \cdot (Z_{1 \; 3'}^2)^{Z_{1
\; 3}^2 Z^{-2}_{2', 3 \; 3'}} \\ & \cdot {(Z_{2 \; 2'})}^{Z^{-2}_{1
\; 1', 2} \bar{Z}^2_{2, 3 \; 3'}} \cdot {(Z^3_{2, 3 \; 3'})}^
{Z^2_{2 \; 2'}} \cdot  Z^3_{1 \; 1', 2'} \cdot Z_{2 \; 2'}.
\nonumber
\end{eqnarray}

However, the resulting factorization $\widetilde{\Delta}$ is not a
braid monodromy factorization of a curve $S' \in V(6,6,4)$, due to
the existence of extra branch points outside $U$. Let $D$ be a disc in
$\C_{x_0}$ containing all the six points in the fiber. Define the
forgetting homomorphisms:
$$
1 \leq i \leq 3\, \, f_i : B_{6} [D, \{1,1',2,2', 3,3' \}]
\rightarrow B_2 [D, \{i,i' \}].
$$

It is clear that if $\widetilde{\Delta}$ were a BMF, then for all
$i,\, \mbox{deg}(f_i(\widetilde{\Delta})) = 2$, by Remark \ref{remArt}.
However, this is not the case in the current situation.  It was
proven in \cite{RobbT} (see also \cite{Robb}), that if
deg$(f_i(\widetilde{\Delta})) = k < 2$, then there are $(2-k)$ extra
branch points, and so there is a contribution of the factorization
$\prod\limits^{2-k}_{m=1} Z_{i,i'}$ to $\widetilde{\Delta}$ (by
contribution we mean that we multiply $\widetilde{\Delta}$ from the
right by these $Z_{i,i'}$'s).

It is easy to see that $\mbox{deg}(f_2(\widetilde{\Delta})) = 2$. In
addition, we have the following

\begin{lemma}
${\emph{deg}}(f_1(\widetilde{\Delta})) =
{\emph{deg}}(f_3(\widetilde{\Delta})) = 1$.
\end{lemma}

\begin{proof}
We prove the lemma only for $f_1$; the proof for $f_3$ is identical.
The braids coming from the nodes are sent by $f_1$ to $Id$, and also the braids
${(Z_{2 \; 2'})}^{Z^{-2}_{1 \; 1', 2} \bar{Z}^2_{2, 3 \; 3'}} $, $ {(Z^3_{2,
3 \; 3'})}^ {Z^2_{2 \; 2'}}$, $ Z_{2 \; 2'}$.
By \cite[Lemma 2, (i)]{19}, we see that deg$(f_1(Z^{3}_{1\,1',2})) = 1$.
\end{proof}

Multiplying $\widetilde{\Delta}$ from the right by $Z_{1,1'} \cdot Z_{3,3'}$ we get a factorization of a curve with
four nodes, six cusps and four branch points, which is isotopic to the branch curve of the Cayley
surface $C$. By \cite[Theorem VI.2.1l]{16} this is indeed the braid monodromy factorization of the branch curve.

\begin{figure}[ht]
\epsfysize=10cm 
\begin{center}
\epsfbox{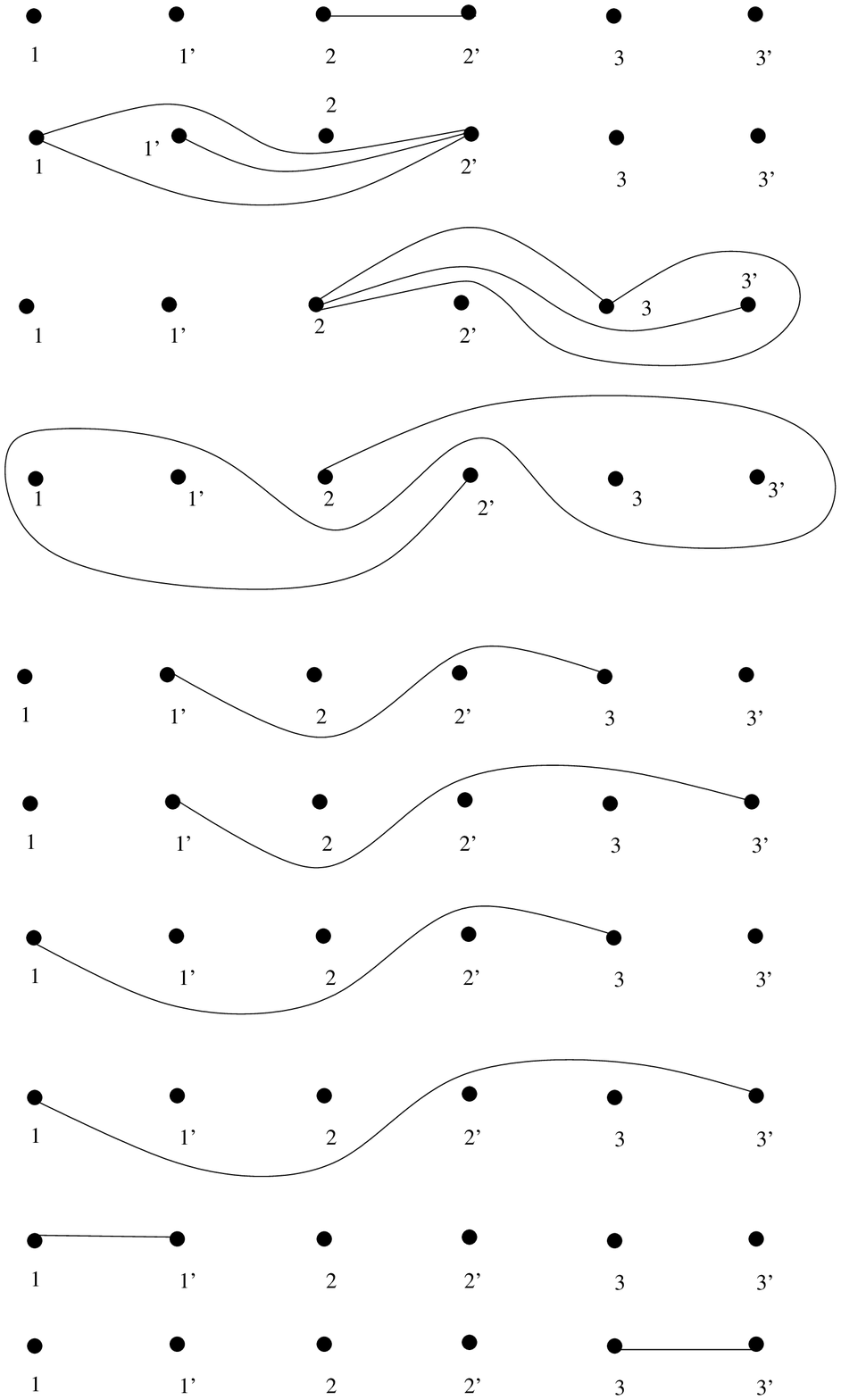}
\end{center}
\caption{}\label{mondel1}
\end{figure}

\begin{theorem}
The braid monodromy factorization of $S$ is given in (\ref{del1})
and its factors are represented by paths in Figure \ref{mondel1}.
\begin{eqnarray} \label{del1}
{\Delta^2_S} = & {(Z_{1' \; 3}^2)}^{Z^{-2}_{2', 3 \; 3'}}
\cdot (Z_{1' \; 3'}^2)^{Z_{1' \; 3}^2 Z^{-2}_{2', 3 \; 3'}} \cdot
{(Z_{1 \; 3}^2)}^{Z^{-2}_{2', 3 \; 3'}} \cdot (Z_{1 \;
3'}^2)^{Z_{1 \; 3}^2 Z^{-2}_{2', 3 \; 3'}} \\ & \cdot {(Z_{2 \;
2'})}^{Z^{-2}_{1 \; 1', 2} \bar{Z}^2_{2, 3 \; 3'}} \cdot {(Z^3_{2,
3 \; 3'})}^ {Z^2_{2 \; 2'}} \cdot  Z^3_{1 \; 1', 2'} \cdot Z_{2 \;
2'} \cdot Z_{1,1'} \cdot Z_{3,3'}. \nonumber
\end{eqnarray}
\end{theorem}

Note that the first, the fourth and the last two paths correspond
to braids induced from the branch points, the second and third ones correspond
to the cusps and the rest correspond to the nodes.

\subsection{The fundamental group $\pi_1(\C^2 - S)$}

The Van Kampen Theorem \cite{vk} states that there is a ``good"
geometric base $\{\gamma_{j}\}$ of $\pi_{1}(\C_{x_{0}} - S \cap
\C_{x_{0}} , *)$ (where $\C_{x_0}$ is the fiber of the projection
$\pi|_{aff}$ above $x_0$), such that the group $\pi_{1}(\C^2 - S,
*)$ is generated by the images of $\{ \gamma_{j} \}$ in
$\pi_{1}(\C^2 - S, *)$ with the following relations:
$\{(\varphi(\delta_{i})) \hspace{0.1cm} (\gamma_{j}) = \gamma_{j}
\hspace{0.1cm} \forall i,j\}$. We recall that
$$
\pi_{1}(\C \P^2 - \bar{S}) \simeq \pi_{1}(\C^2 - S) / \langle
\prod{\gamma_{j}}\rangle.
$$

\begin{notation}
$[x,y] \doteq xyx^{-1}y^{-1},\,\,\langle x,y \rangle \doteq xyxy^{-1}x^{-1}y^{-1}.$
\end{notation}

Recall that $S$ has only branch points, nodes and cusps (when the
cusp is locally defined by the equation $y^2 = x^3$). Denote by $a$
and $b$ the two points of $S$ at a neighborhood of a singular point
on the fiber $\C_{x_0}$. Let $\gamma_{a}, \gamma_{b}$ be two
non-intersecting loops in $\pi_{1}(\C_{x_{0}} - S \cap \C_{x_{0}} ,
*)$ around the intersection points $a$ and $b$.
Then by the van Kampen Theorem, we have the relation $\langle
\gamma_{a},\gamma_{b} \rangle =  1$ induced from a cusp, the
relation $[\gamma_{a},\gamma_{b}] =  1$ induced from a node and the
relation $\gamma_a = \gamma_b$ induced from a branch point.

\begin{theorem}
The fundamental group $\pi_1(\C^2 - S)$ is generated by
$\gamma_1, \gamma_{2}, \gamma_{3}$ subject to the relations%
\begin{eqnarray}
{}\langle\gamma_{1}, \gamma_{2}\rangle & = & e \label{fin1}\\
{}\langle\gamma_{2}, \gamma_{3}\rangle & = & e \label{fin2}\\
{}[ \gamma_{2}, \gamma_{1}^2 \gamma_{3}^2] & = & e \label{fin3}\\
{}[ \gamma_{1}, \gamma_{2}^{-1} \gamma_{3} \gamma_2] & = & e.
\label{fin4}
\end{eqnarray}

The group $\pi_1(\C\P^2 - \bar{S})$ has relations (\ref{fin1}),
(\ref{fin2}), (\ref{fin4}) and an additional relation
\begin{eqnarray}\label{relpro}
\gamma_{3}^2 \gamma_{2}^2 \gamma_{1}^2 & = & e.
\end{eqnarray}
\end{theorem}

\begin{proof}
By the above explanation and Figure \ref{mondel1}, we have the
following relations
\begin{eqnarray}
{}\gamma_2 & = & \gamma_{2'} \label{22'} \\
{}\langle\gamma_{1}, \gamma_{2'}\rangle = \langle \gamma_{1'},
\gamma_{2'}\rangle =
  \langle\gamma_{1'} \gamma_1 \gamma_{1'}^{-1}, \gamma_{2'}\rangle &=& e
\label{trip12} \\
{}\langle\gamma_{2'} \gamma_{2} \gamma_{2'}^{-1}, \gamma_3\rangle
= \langle\gamma_{2'} \gamma_{2} \gamma_{2'}^{-1},
\gamma_{3'}\rangle =
  \langle\gamma_{2'} \gamma_{2} \gamma_{2'}^{-1}, \gamma_{3'} \gamma_3
\gamma_{3'}^{-1} \rangle & = & e \label{trip23} \\
{}\gamma_{1}^{-1} \gamma_{1'}^{-1} \gamma_{2'}^{-1} \gamma_{3'}
\gamma_{3} \gamma_{2'} \gamma_{2} \gamma_{2'}^{-1} \gamma_3^{-1}
\gamma_{3'}^{-1} \gamma_{2'}  \gamma_{1'}
\gamma_{1} & = & \gamma_{2'} \label{bpoint} \\
{}[\gamma_{1}, \gamma_{2'}^{-1} \gamma_3 \gamma_{2'}] =
[\gamma_{1'}, \gamma_{2'}^{-1} \gamma_3^{-1} \gamma_{3'}
\gamma_{3} \gamma_{2'}] & = & e \label{13a} \\
{}[\gamma_{1}, \gamma_{2'}^{-1} \gamma_{3} \gamma_{2'}] =
[\gamma_{1},\gamma_{2'}^{-1} \gamma_3^{-1} \gamma_{3'} \gamma_{3}
\gamma_{2' }] &
= & e \label{13b} \\
{}\gamma_{1} & = & \gamma_{1'} \label{11'} \\
{}\gamma_{3} & = & \gamma_{3'}. \label{33'}
\end{eqnarray}

We want to simplify this presentation. By (\ref{22'}) and
(\ref{11'}), relation (\ref{trip12}) gets the form (\ref{fin1}).
By (\ref{22'}) and (\ref{33'}), relation (\ref{trip23}) gets the
form (\ref{fin2}). Relation  (\ref{bpoint}) is transformed (by
(\ref{fin1}), (\ref{22'}), (\ref{11'}) and (\ref{33'})) to
(\ref{fin3}), and relations (\ref{13a}) and (\ref{13b}) are
transformed (by (\ref{22'}), (\ref{11'}) and (\ref{33'})) to
(\ref{fin4}).

In order to get the group $\pi_1(\C\P^2 - \bar{S})$, we add the
projective relation $\gamma_{3'} \gamma_{3} \gamma_{2'} \gamma_{2}
\gamma_{1'} \gamma_{1} = e$, which is transformed to $\gamma_{3}^2
\gamma_{2}^2 \gamma_{1}^2 = e$. Therefore relation (\ref{fin3}) is
omitted and $\pi_1(\C\P^2 - \bar{S})$ is generated by $\gamma_1,
\gamma_{2}, \gamma_{3}$ with  relations (\ref{fin1}),
(\ref{fin2}), (\ref{fin4}) and (\ref{relpro}).

We note that if we consider the relations which are derived from
the complex conjugates of the braids of Figure \ref{mondel1}, we
gain no new relations, therefore the presentation is complete.
\end{proof}

\begin{remark}
\emph{Since we deal with a singular cubic surface, we  note that a
result of Zariski \cite{Z} for a smooth cubic surface in $\C\P^3$
was generalized by Moishezon \cite{Mo} for any degree. Let
$\bar{S}_n$ (resp. $S_n$) be the branch curve of a smooth surface of degree $n$
in $\C\P^3$ in $\C\P^2$ (resp. $\C^2$). Moishezon proved that
$$
\pi_{1}(\C^2 - S_n) \isom B_n  \ \ \ \mbox{and} \ \ \
\pi_{1}(\C\P^2 - \bar{S}_n) \isom B_n/Center(B_n)
$$}
\end{remark}

\section{Finding the fundamental group $\pi_1(C - Sing)$}\label{sec:6}

In this section we give two different ways to find the fundamental
group of the complement of the singular points in the Cayley
surface.

\subsection{Using a lifting to the Mapping Class Group}\label{5.1}

The projection $\pi$ defines a pencil of lines on $\C\P^2$.
Considering the preimages of these lines under the projection of $C$
onto $\C\P^2$, we obtain a pencil of elliptic curves on $C$,
intersecting transversely at the base locus, namely, three smooth
points (the preimages in $C$ of the pole of the projection $\pi$).
This pencil has eight nodal fibers, of which four pass through the
singular points of $C$ and the four others pass through the
preimages of the branch points of $S$ with respect to the projection
$\pi$. The monodromy of this fibration can be encoded by a
factorization in a mapping class group, which can be obtained from
the braid monodromy of $S$ by a simple lifting algorithm. See
Section 5.2 of \cite{AK} and Section 3.3 of \cite{AuG} for details.

Among the various factors of the factorization (\ref{del1}), those corresponding to
cusps of $S$ (i.e. $Z^3_{1 \; 1', 2'}$ and ${(Z^3_{2, 3 \; 3'})}^
{Z^2_{2 \; 2'}}$) lie in the kernel of the lifting homomorphisms
and do not contribute to the monodromy of the elliptic pencil.
This is because the preimage of the fiber of $\pi$ through a cusp
of $S$ is actually a smooth elliptic curve.

To determine the liftings of the other factors, we view the fiber
$E$ of the elliptic pencil as a triple cover of a line in $\C\P^2$
(the reference fiber of $\pi$ on which the braid monodromy acts)
branched at six points (the points where $S$ intersects the
considered line), which we label $1,1',2,2',3,3'$ as before. Each of
these branch points corresponds to a simple ramification, i.e.,
involving only two of the three sheets of the covering. We label
these sheets by elements of $\{1,2,3\}$. We need to find the
monodromy epimorphism $\pi_1(\C\P^2 - \bar{S}) \rightarrow Sym_3$
such that if $\langle a,b \rangle = 1$ for a pair of generators $a,b
\in \pi_1(\C\P^2 - \bar{S})$, then $a$ and $b$ would be sent to two
non-commuting transposition $(i\, j)$ and $(j\, k)$ (see, for
example, \cite{Mo0} for the explicit conditions imposed on the
monodromy epimorphism). The only epimorphism, up to renumeration of
the sheets, is the one that maps $\gamma_1$ and $\gamma_{1'}$ to the
transposition $(23)$, $\gamma_2$ and $\gamma_{2'}$ to $(13)$, and
$\gamma_3$ and $\gamma_{3'}$ to $(12)$.
%
%

 The lifting homomorphism (see \cite{AK,AuG}) maps
the halftwist $Z_{1\,1'}$ to a positive Dehn twist along the simple
closed loop on $E$ formed by the two lifts of the supporting arc of
the halftwist in sheets $2$ and $3$ of the covering. This can be
done similarly for the other halftwists appearing in the braid
monodromy factorization. Because the halftwists $Z_{i\,i'}$ ($1 \le
i \le 3$) have disjoint supporting arcs, the corresponding Dehn
twists $\tau_{i\,i'}$ also have disjoint supporting loops; moreover
it is easy to check that these loops are homotopically non-trivial
(see below). Since disjoint non-trivial simple closed loops on an
elliptic curve are homotopic, the Dehn twists $\tau_{i\,i'}$
correspond to mutually homotopic vanishing cycles, and represent the
same element in the mapping class group $\mathrm{Map}_1=SL(2,\Z)$.
We call $\alpha\in \pi_1(E)$ the homotopy class of these vanishing
cycles, and $\tau_\alpha$ the corresponding Dehn twist. Next we
observe that the support of the halftwist
$t={(Z_{2\,2'})}^{Z^{-2}_{1\,1',2}\bar{Z}^2_{2,3\,3'}}$ (Figure
\ref{mondel1}, fourth line) is also disjoint from those of
$Z_{1\,1'}$ and $Z_{3\,3'}$, which indicates that the corresponding
vanishing cycle again represents the homotopy class $\alpha$ in
$\pi_1(E)$.

For a generic projection of a smooth surface, the nodes of the
branch curve correspond to smooth fibers of the pencil, and the
corresponding braid monodromies lie in the kernel of the lifting
homomorphism. However, in our case the four nodes of the branch
curve correspond to nodal fibers of the pencil; the corresponding
braid monodromies are squares of liftable halftwists, which lift
to the squares of the corresponding Dehn twists.

We first consider $\nu=(Z_{1'\,3})^{Z^{-2}_{2',3\,3'}}$ (Figure
\ref{mondel1}, fifth line): the supporting arc of this halftwist
intersects that of $Z_{1\,1'}$ only once, at the common endpoint
$1'$. Hence, the double lifts of the supporting arcs (which give
the supporting loops of the corresponding Dehn twists) intersect
transversely exactly once (at the branch point which lies above
$1'$). Calling $\beta\in \pi_1(E)$ the homotopy class of the
vanishing cycle corresponding to the lift of
$(Z_{1'\,3})^{Z^{-2}_{2',3\,3'}}$, the intersection number
$\alpha\cdot \beta=1$ implies that $\alpha$ and $\beta$ form a
basis of $\pi_1(E)\simeq \Z^2$ (and confirms that the vanishing
cycles are indeed not homotopically trivial as claimed above). The
same argument could have been used considering $Z_{3\,3'}$ (whose
support intersects that of $\nu$ once at the common end point
$3$), or $Z_{2\,2'}$ or $t$ instead of $Z_{1\,1'}$ (in that case
the supporting arcs intersect transversely once at interior
points, but their double lifts each live in only two of the three
sheets of the covering, and it is easy to check that the
supporting loops of the corresponding Dehn twists intersect
transversely once). In any case, we conclude that the braid
monodromy factor $\nu^2$ lifts to $\tau_\beta^2$, the square of
the positive Dehn twist about a loop in the homotopy class
$\beta$.

Finally, the three other nodes of $S$ ($(Z_{1' \; 3'}^2)^{Z_{1' \;
3}^2 Z^{-2}_{2', 3 \; 3'}}, {(Z_{1 \; 3}^2)}^{Z^{-2}_{2', 3 \; 3'}}$
and $(Z_{1 \; 3'}^2)^{Z_{1 \; 3}^2 Z^{-2}_{2', 3 \; 3'}}$)
correspond to the conjugates of $\nu^2$ by the braids
$Z_{3\,3'}^{-1}$, $Z_{1\,1'}^{-1}$, and
$Z_{1\,1'}^{-1}Z_{3\,3'}^{-1}$, respectively (Figure \ref{mondel1},
sixth, seventh, eighth lines). Applying the lifting homomorphism, we
obtain that the corresponding mapping class group elements are
respectively
$(\tau_\beta^2)^{\tau_\alpha^{-1}}=\tau_{\beta-\alpha}^2$,
$(\tau_\beta^2)^{\tau_\alpha^{-1}}=\tau_{\beta-\alpha}^2$, and
$(\tau_\beta^2)^{\tau_\alpha^{-2}}=\tau_{\beta-2\alpha}^2$. In other
words, the vanishing cycles represent respectively the homotopy
classes $\beta-\alpha$, $\beta-\alpha$, and $\beta-2\alpha$.

In conclusion, the mapping class group monodromy factorization of
the elliptic pencil (in $\mathrm{Map}_1$) is
\begin{equation}\label{delfin}
\mathrm{Id}=\tau_\alpha\cdot \tau_\alpha\cdot \tau_\beta^2\cdot
\tau_{\beta-\alpha}^2 \cdot \tau_{\beta-\alpha}^2 \cdot
\tau_{\beta-2\alpha}^2\cdot \tau_\alpha\cdot \tau_\alpha.
\end{equation}

As a verification, one can consider the isomorphism
$\mathrm{Map}_1\simeq SL(2,\Z)$ given by the action on $H_1(E,\Z)$,
working i the basis $\{\alpha,\beta\}$. Then, recalling that the
action of a Dehn twist on homology is given by
$[\tau_\delta(\gamma)]=[\gamma]+ ([\delta]\cdot[\gamma])[\delta]$,
we have
$$\tau_\alpha=\begin{pmatrix} 1&1\\0&1 \end{pmatrix},\quad
\tau_\beta=\begin{pmatrix} 1&0\\-1&1 \end{pmatrix},\quad
\tau_{\beta-\alpha}=\begin{pmatrix} 2&1\\-1&0 \end{pmatrix},\quad
\tau_{\beta-2\alpha}=\begin{pmatrix} 3&4\\-1&-1 \end{pmatrix},$$ and
the identity (\ref{delfin}) indeed holds (recalling that our
products are written in the braid order, i.e., with composition from
left to right,
while the usual product of matrices is a composition from right to left).\\

Let us recall the quasi-projective Lefschetz Hyperplane Section
Theorem: {\it Let $X:=Y- Z,$ (of dimension $d$) where $Y$ is an
algebraic subset of the complex projective space $\C\P^n, \, n\geq
2,$ and where $Z$ is an algebraic subset of $Y.$ Let $L$ be a
projective hyperplane which is in generic position with respect to
$X.$ If $X$ is nonsingular, then the natural maps $$H_q(L\cap X)\to
H_q(X) \quad {\rm and} \quad \pi_q(L\cap X,\ast )\to \pi_q(X)$$ are
bijective for $0\leq q \leq d-2$ and surjective for $d-1$ (see
\cite{HL} and \cite{GM})}. In our case $X$ is the Cayley cubic $C$
minus the singular locus. The above generically chosen central
projection $\pi:\C\P^2- P_0\to \C\P^1$ (where $P_0$ is the center of
the projection) defines a pencil of lines in $\C\P^2$ which lifts to
a generic pencil of planes in $\C\P^3$ whose axis is a line $M$. The
above curve $E$ is then the intersection of $C$ with a generic
member $L$ of the pencil of planes. It follows that the natural map
$$\pi_1((C-Sing) \cap L)=\pi_1(E)\to \pi_1(C-Sing)$$ is a
surjection. In particular, the fundamental group of $C$ is abelian,
so we can work with homology groups instead of fundamental groups.
It follows that we are left to determine the kernel of the natural
map $$ \varphi:H_1(E)\to \pi_1(C-Sing).$$

Let $L_i$ be the  exceptional hyperplanes of the above pencil (the
planes for which $L_i\cap (C-Sing)$ are not isotopic to $E$). In
\cite{Cheniot}, Cheniot  defines homological variation operators
$${\rm var}_{i,q}: H_q(X\cap L, M\cap X) \to  H_q(X\cap L),\quad i
\in I,$$ by patching each relative cycle on $X\cap L$ modulo $M\cap
X$ with its transform by monodromy around the exceptional lines. It
is then shown in \cite{Cheniot}, that
$$ {\rm Kernel}(\varphi)=\sum_{i\in I} {\rm Im}({\rm var}_{i,q}).$$
If we choose a basis $\{\alpha,\beta\}$ of $H_1(E)\simeq H_1(E,M\cap
X)$ as above, then the image of ${\rm var}_{i,q}$ is nothing else
but  the image of the lifted braid, viewed as an element in the
mapping class group ${\rm Map}_1.$ (This can be seen by unravelling
the definitions of Cheniot and the above construction of the mapping
class group factorization, see also \cite{CheniotEyral}, Section 2).
It follows that in our case $\pi_1(X)=\pi_1(C-Sing)$ is isomorphic
to the quotient of $\pi_1(E)=\Z\alpha\oplus \Z\beta$ by the
relations $\gamma= \rho_*(\gamma)$ for every $\gamma\in \pi_1(E)$
and for every factor $\rho$ in the mapping class group factorization
\eqref{delfin}.
 The relation
$\tau_\alpha(\beta)=\beta$ implies that $[\alpha]$ is trivial in
$\pi_1(C-Sing)$, while the relation $\tau_\beta^2(\alpha)=\alpha$
implies that $2[\beta]$ is trivial; hence, $\pi_1(C-Sing)$ is a
quotient of $\Z/2$. It follows from the monodromy factorization
given in Formula \eqref{delfin} that
 for every element $\rho$ of the monodromy subgroup, the image of
$\rho-\mathrm{Id}$ is in the span of $\alpha$ and $2\beta$.
Therefore,  we get no further relations, and we recover that
$\pi_1(C-Sing)=\Z/2$.

\subsection{Using the RMS method}\label{5.2}

We find $\pi_1(C - Sing)$ by a second method, using the
Reidemeister-Schreier algorithm. We follow the method proposed in
\cite{Manf}. We recall this method briefly.

Denote by $Gr^*_{d,n}$ the set of the graphs with $d$ labelled
vertices and $n$ labelled edges. Assume we have a homomorphism $f:
\pi_1(\C^2 - S) \ri Sym_d$ and let $\g_1,...,\g_n$ be generators of
$\pi_1(\C^2 - S)$. Denote by $g$ the homomorphism from the free
group with $n$ generators $F_n = \langle \bar{\g}_1,...,\bar{\g}_n
\rangle$ to $\pi_1(\C^2 - S)$ such that $g(\bar{\g}_i) = \g_i,\,
\forall\,i = 1,...,n$.

Assume we have a homomorphism $\bar{f}: F_n \ri Sym_d$ such that
$\bar{f}(\bar{\g}_i)$ is a transposition $\forall\,i = 1,...,n$.
So we can associate to it a graph $\G_{\bar{f}} \in Gr^*_{d,n}$ in
the following way. If $\bar{f}(\bar{\g}_i) = (h,k)$, then the edge
$i$ will have the vertices $h$ and $k$. Given a monodromy map
$f:\pi_1(\C^2 - S) \ri Sym_d$, the monodromy graph associated to
it is the graph $\G = \G_{\bar{f}}$, where $\bar{f}$ is the
lifting of $f$ to $F_n$ under the map $g$.

In order to compute $\pi_1(C - Sing)$, let us consider the
projection $\pi|_{C - D}: C - D \ri \C^2 - S$, where $D =
\pi^{-1}(S)$ and $D = 2R + F$ ($R$ is the ramification locus of
$\pi$). As this is an unramified cover, we can identify $\pi_1(C -
D)$ with the subgroup of $\pi_1(\C^2 - S)$ given by those elements
$\g$ such that $\g$ stabilizes a vertex of $\G$ (i.e.,  the $\g$'s
such that $f(\g)(j)=j$ for a fixed vertex $j$ of $\G$).

Explicitly, taking a base point in $C - D$ to be the preimage of
the base point of $\C^2 - S$ lying in the sheet labelled 1, then a
loop $g$ in $\C^2 - S$ lifts to an arc in $C - D$, whose other end
point is the preimage in the sheet $f(g)(1)$ ; hence we obtain a
closed loop in $C - D$ if and only if $f(g)$ maps 1 to 1.

Let us fix a numeration on $\G$, and let $\G'$ be a maximal subtree.
By abuse of notation, let us denote by $\g_i$ the edges of $\G$.
\begin{definition}
A sequence $c = (k_j)_{j = 1,..,l}$ of distinct edges of $\G$ such
that the edge $k_i$ intersects the edge $k_{i+1}$ only in a single
vertex is called a chain of $\G$ (of length $l$). A $p$-chain is a
chain with $p$ as a starting vertex. A $p,q$-chain is a chain with
$p$ as a starting vertex and $q$ as an ending vertex.
\end{definition}

If $c = (k_j)_{j = 1,..,l}$ is a 1-chain in $\G'$, then set $\g_c
= \g_{k_1}...\g_{k_l}$, and if $c$ is the trivial 1-chain, then
set $\g_c = id$. The set of all $\sigma \in Sym_d$ such that
$\sigma = f(\g_c)$ for $c$ a 1-chain in $\G'$ is a complete set of
representatives for left cosets of the stabilizer of the vertex 1
in $Sym_d$ :  if $c$ is a $1,q$-chain in $\G'$ then $f(\g_c)(1) =
q$.

So in order to calculate $\pi_1(C - D)$, we apply the
Reidemeister-Schreier method (\cite{MKS}) to the Schreier set $RS =
\{\g_c |\, c\,\, \mbox{is a 1-chain in}\,\, \G'\}$ to get the
following proposition \cite[Prop. 6.1]{Manf}:

\begin{proposition}
$\pi_1(C - D)$ is generated by $\eta_{c,k} =
\g_c\g_k(\overline{\g_c\g_k})^{-1}$, where $c$ is a $1$-chain in
$\G'$, $k$ an edge of $\G$ (such that $c \cup \{k\}$ is not a
$1$-chain in $\G'$), and is defined by the relators $\g_cR\g_c^{-1}$
(written in terms of $\eta$'s), where $c$ is a $1$-chain in $\G'$,
$R$ is a relator of $\pi_1(\C^2 - S)$.
\end{proposition}

Reading the proof from \cite{Manf}, we see that in order to obtain a
presentation for $\pi_1(C - Sing)$, we must quotient by the normal
subgroup generated by all loops around the components $D = 2R+F$.
The loops around the components of $R$ are those $\eta_{c,k}$ with
$k$ equal to the last edge of $c$ and loops $\eta_{c,k}\eta_{c',k}$
in case $k$ is an edge of $\G$ such that $c \cup \{k\}$ is not a
1-chain in $\G'$, and $\g_{c'} = \overline{\g_c\g_k} \in RS$; while
the loops around the components of $F$ are those $\eta_{c,k}$ with
$k$ an edge which does not pass through the ending vertex of $c$.

\begin{figure}

\begin{minipage}{\textwidth}

\begin{center}

\epsfbox{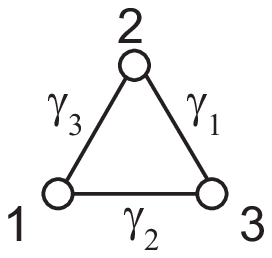}

\end{center}

\end{minipage}

\caption{The graph $\G$}\label{mich1}

\end{figure}

Using this proposition, we can find $\pi_1(C - Sing)$ in our case.
Recall that the monodromy maps $\g_1$ to the transposition
$(2,3)$, $\g_2$ to $(1,3)$ and $\g_3$ to $(1,2)$. In this way, we
create the map $f:\pi_1(\C^2 - S) \ri Sym_3$, and thus we can
associate to it the graph $\G = \G_{\bar{f}}$ in Figure
\ref{mich1}.

Denote by $\G'$ the maximal subtree composed of the edges $\g_1$
and $\g_3$ (and the vertices $\{1,2,3\}$). Let RS =
$\{Id,\g_3,\g_3\g_1\}$ be the Schreier set of the stabilizer of
the vertex 1. Let us denote $\g_{c_1} = Id,\g_{c_2} = \g_3,
\g_{c_3} = \g_3\g_1$ , $k_i = i$, and let $\eta_{i,j} =
\g_{c_i}\g_{k_j}(\overline{\g_{c_i}\g_{k_j}})^{-1}$ be the
generators of $\pi_1(C - D)$.

We get the following generators:\\
$\eta_{1,1} = Id\cdot\g_1(\overline{Id\cdot\g_1})^{-1} = \g_1$
\newline $\eta_{1,2} = Id\cdot\g_2(\overline{Id\cdot\g_2})^{-1} =
\g_2\g_1^{-1}\g_3^{-1}$ \newline $\eta_{1,3} =
Id\cdot\g_3(\overline{Id\cdot\g_3})^{-1} = \g_3\g_3^{-1} = Id$
\newline
$\eta_{2,1} = \g_3\cdot\g_1(\overline{\g_3\cdot\g_1})^{-1} =
\g_3\g_1\g_1^{-1}\g_3^{-1} = Id$ \newline $\eta_{2,2} =
\g_3\cdot\g_2(\overline{\g_3\cdot\g_2})^{-1} =
\g_3\g_2\g_3^{-1}$\newline $\eta_{2,3} =
\g_3\cdot\g_3(\overline{\g_3\cdot\g_3})^{-1} = \g_3^2$\newline
$\eta_{3,1} = \g_3\g_1\cdot\g_1(\overline{\g_3\g_1\cdot\g_1})^{-1}
= \g_3\g_1^2\g_3^{-1}$\newline $\eta_{3,2} =
\g_3\g_1\cdot\g_2(\overline{\g_3\g_1\cdot\g_2})^{-1} =
\g_3\g_1\g_2$\newline $\eta_{3,3} =
\g_3\g_1\cdot\g_3(\overline{\g_3\g_1\cdot\g_3})^{-1} =
\g_3\g_1\g_3\g_1^{-1}\g_3^{-1}$.\newline

Thus we have 7 generators; to find $\pi_1(C - Sing)$ we have to
eliminate all loops around $D$. That is, we quotient by the
generators $\eta_{1,1}, \eta_{2,2}, \eta_{2,3}, \eta_{3,1},
\eta_{3,3}$, leaving us with $\eta_{1,2}$ and $\eta_{3,2}$. However,
$\eta_{1,2}\cdot\eta_{3,2} = \g_2^2$, which also represents a loop
around $D$, and thus $\g_2^2 =Id$. Thus, $\pi_1(C - Sing)$ is
generated by one element $\g_2$,  and isomorphic to
$\mathbb{Z}/2\Z$.

\section{Acknowledgements}

This work was initiated while the first named author was hosted at
the Mathematics Institute, Erlangen-N\"urnberg University, Germany.
Thanks are given to the Institute and to the host Wolf Barth, who
presented the problem and gave motivation for this work. We are
grateful to Ron Livne and Yann Sepulcre for fruitful discussions.
Thanks are also given to Denis Auroux for helping us complete
Section \ref{5.1}. We also  thank Anatoly Libgober and Sandro
Manfredini for helpful discussions in Section \ref{5.2}.

\end{document}